\documentclass[12pt,english]{article}
\usepackage{amsthm,amsmath,natbib}
\usepackage{dsfont}
\usepackage{graphicx}
\usepackage{amsmath,amssymb,amsfonts,hhline}
\usepackage[hmargin=3cm,vmargin=3cm]{geometry}
\usepackage{epsfig}
\usepackage{epstopdf}
\newtheorem{example}{Example}
\newtheorem{definition}{Definition}
\newtheorem{remark}{Remark}
\newcommand{\sign}{\hbox{sign}}

\begin{document}
\title{Inference on the marginal distribution of clustered data with informative cluster size}
\author{J. Nevalainen \\ University of Turku, Finland \\ email: {jaakko.nevalainen@utu.fi} \\ \\
S. Datta \\ University of Louisville, USA \\email: {somnath.datta@louisville.edu}
\\ \\
H. Oja \\  University of Tampere, Finland \\email: {hannu.oja@uta.fi}}

\maketitle

\begin{abstract}
In spite of recent contributions to the literature, informative  cluster size settings
are not well known and understood. In this paper,
we give a formal definition of the problem and describe it from different viewpoints.
Data generating mechanisms, parametric and nonparametric models are considered in light of examples.
Our emphasis
is on nonparametric and robust approaches to the inference on the marginal distribution.
Descriptive statistics and parameters of interest are defined as functionals and they are accompanied with a generally
applicable testing procedure. The theory is illustrated with an example on patients with incomplete spinal cord
injuries.
\end{abstract}


\section{Introduction}

Clustered data problems are encountered everywhere in biomedical
research and, not surprisingly, the statistical methods involving
the analysis of cluster correlated data have been subject to
intensive research even until today. A typical
situation is that instead of sampling $N$ independent and
identically distributed (i.i.d.) random variables, the researcher
samples observations in (say) $M$ clusters with known cluster
memberships. Observations within a cluster tend to be similar in some
way but can be assumed independent across clusters.

To be more specific, write  $Y_{i1},...,Y_{iN_i}$ for the $N_i$
observations in the $i$th cluster, $i=1,...,M$. Let $X_{ij}$ be a
possible vector of (random or fixed) explanatory variables for the
response value $Y_{ij}$, $i=1,...,M$; $j=1,...,N_i$. The cluster
sizes $N_1,\ldots,N_M$  are often simply thought to be fixed design
constants. In the linear regression  model it is then assumed
that, for a correct value $\beta$,
\begin{itemize}
\item[(i)] the marginal distributions of
$\epsilon_{ij}=Y_{ij}-\beta'X_{ij}$ are all the same,
\item[(ii)] all the
bivariate distributions of $(\epsilon_{ij},\epsilon_{ij'})$, $j\ne
j'$ are the same, and
\item[(iii)] $\epsilon_{ij}$ and $\epsilon_{i'j'}$,
$i\ne i'$ are independent.
\end{itemize}
If multivariate normality
of the random errors $\epsilon_{ij}$ can be assumed, for example, the most popular
technique for valid statistical inference for the parameter $\beta$
is to employ mixed models with cluster effects as random effects.
Alternatively, one can work out the variance terms for different
test statistics, and modify the tests accordingly. Introduction of
weights (indirectly present in likelihood inference for mixed
models) can potentially improve the efficiency of the analysis, but
variance adjusted test and estimating procedures based on different
weighting schemes all provide valid statistical inference in this
model. If  $N_1,...,N_M$ are random  and the joint
distribution of the random errors $\epsilon_{ij}$ does not depend on
$N_1,...,N_M$, it is still reasonable to assume
that (i)--(iii) hold.

A much more complex  setting arises when the cluster size may have
an influence on the measured values, or vice versa, or possibly they
are both influenced by a third, unobservable latent variable. The setting
is termed \emph{informative cluster size}, because the cluster
size---which is now also a random variable---could then carry
information about the quantities or parameters of interest. Recent examples
of informative cluster size problems in the biostatistical literature include:
\begin{itemize}
\item volume-outcome studies \citep{Panageas2007} where specialized
surgeons treating many patients
may have better outcomes than those treating few patients;
\item periodontal
studies \citep{WilliamsonDattaSatten2003,Wang2011} where patients with fewer teeth tend to have poorer
condition of the still remaining teeth;
\item radiation
toxicity studies \citep{Datta2008}, where the number of
measurements on successive measurement on an individual depends on the number of
radiation therapies, which in turn depend on the
underlying severity of cancer.
\end{itemize}
More examples are provided in \citet{Dunson2003}, \citet{Williamson2007} and \citet{Williamson2008}, among others.

\citet{Hoffman2001}, in their original paper, defined \emph{nonignorable cluster size}
as any violation of the property that $E(Y_{ij}|X_{ij},N_i)=E(Y_{ij}|X_{ij})$ in the framework of generalized estimation equations.
As the recent interest in informative cluster size problems has gone far beyond that particular setting, it makes sense to define
the concept in a more general way.
\begin{definition}
We say that cluster size is noninformative if
$$P(Y_{ij}\leq y|N_i=k)=P(Y_{ij}\leq y),\ \ k=1,2,...; j=1,...,k$$
Otherwise, it is called informative.
 Cluster size is  conditionally noninformative if
$$P(Y_{ij}\leq y|X_{ij},N_i=k)=P(Y_{ij}\leq y|X_{ij}),\ \ k=1,2,...; j=1,...,k$$
Otherwise, it is conditionally informative.
\end{definition}
The first part of the definition shows that, in the case of noninformative cluster
size $N_i$ and exchangeable $Y_{i1},...,Y_{iN_i}$,
the characteristics of their common marginal distribution can be estimated consistently in the usual
way, given that the variance terms are corrected appropriately for clustering.
By exchangeability we mean that
$$Y_{ip_1},\ldots,Y_{ip_{N_i}}\sim Y_{i1},\ldots,Y_{iN_i}$$ for all
permutations $(p_1,...,p_{N_i})$ of $(1,...,N_i)$.
Similarly, if the cluster size is conditionally noninformative,
the relationship between $Y_{ij}$ and $X_{ij}$ is not influenced by
cluster size and standard clustered data methods can be used.
The standard approaches are not sufficient if the condition is violated.

The outline of this paper is as follows.  We first
discuss possible data generating mechanisms and appropriate models. They
are illustrated with examples from the literature. In section \ref{sec:functionals}, we
formulate extensions of common quantities of interest---such as quantities of location, scale,
and correlation, and regression coefficients---as alternative functionals, and
discuss appropriate choices among their sample counterparts. Section \ref{sec:testing}
gives a general recipe for constructing tests on the functionals. The theory is illustrated with an example on patients with incomplete spinal cord
injuries in Section \ref{sec:dataexample} and Section \ref{sec:concluding} provides concluding remarks.

\section{Models for informative cluster size problems}

A reasonable basis for statistical modeling is to assume that the measurements
in the $i$th cluster are: the cluster size  $N_i$ and a realization
of a stochastic process
$(Y_{ij})_{j=1}^\infty$, that is,
$ (N_i; Y_{i1},Y_{i2},\ldots ).$
The stochastic process may be finite or infinite, and possibly multivariate.
In practice we only observe  $N_i$ and $Y_i=(Y_{i1},\ldots,Y_{iN_i})'$,
and the observed data consists of
$$V_i=(N_i; Y_{i1},\ldots, Y_{iN_i}),\ i=1,\ldots,M.$$
The cluster variables $V_1,\ldots,V_M$ are assumed to be
independent and identically distributed with some probability
measure $\mathds{P}(V)$.
The goal is  to make valid and efficient statistical inference on the marginal distribution of
$Y_{ij}$'s, which can be a result of an unconventional data generating mechanism.
Any violation of the first condition in Definition 1 means that the cluster size cannot be ignored while
making inference on the marginal distribution of $Y_{ij}$'s, but  needs to be accounted for.

\begin{remark}
A formal way of definition of probability distributions of
$V=(N,Y_1,\ldots,V_N)$, where $N$ is
the cluster size, can be written as follows.

Let $(\Omega,\mathcal{F},P)\,$ be the common probability space on which
$N,Y_1,\,Y_2,\,\cdots$ \ are defined. Then$\,V\,$ is a random element on
$(\Omega,\mathcal{F},P)$ taking values in
$\,\,\mathcal{X}=\,\,\underset{k\geq 1}{\cup}\{k\}\times\Re^k.\,\,\,$The
appropriate $\sigma$-algebra $\,$on $\mathcal{X}$ is the $\sigma$-algebra
generated by the $\pi$-class$\,\,\Pi=$ $\emptyset\cup\{\,\,\{k\}\times
B_1\times\cdots\times B_k:B_i\,=[c_i,d_i)\,$is a semi-open interval in
$\Re;\,1\leq i\leq k,\,k\geq 1\}.\,\,\,$

We define $\mathds{P}_V\,$on such sets by$\,\,\mathds{P}_V(\emptyset)=0\,$ and
\begin{equation*}
\mathds{P}_V(\,\,\{k\}\times B_1\times\cdots\times B_k)=Pr(N=k)Pr(Y_1\in
B_1,\,\cdots\,,\,Y_k\in B_k|N=k).
\end{equation*}
Note that $\mathds{P}_V$ is $\sigma$-additive on $\,\Pi$. Let
$A_i=\{k_i\}\times B_{i1}\times\cdots\times B_{ik_i}\in\Pi,\,i\geq 1,\,$ such
that
$A_i\cap A_j=\emptyset,\,$for $i\neq j,\,$and
$\overset{\infty}{\underset{i=1}{\cup}}A_i=A=\{k\}\times
B_1\times\cdots\times B_k\in\Pi.\,$ Then we must have
$k_i=k,\,\forall i$ and $C=\overset{\infty}{\underset{i=1}{\cup}}C_i$, where
$C=B_1\times\cdots\times B_k,\,C_i=B_{i1}\times\cdots\times
B_{ik};\,$furthermore, $C_i\cap C_j=\emptyset,\,$ for $i\neq j$.

Since the distribution of $(Y_1,\,\cdots,\,Y_k)\,$given $N=k,\,$ denoted
$\mathds{P}^k$, say, is a proper probability distribution,
\begin{equation*}
\mathds{P}^k(\overset{\infty}{\underset{i=1}{\cup}}C_i)=\sum_{i=1}^\infty%
\mathds{P}^k(C_i),
\end{equation*}
and hence
\begin{equation*}
\mathds{P}_V(\,\overset{\infty}{\underset{i=1}{\cup}}A_i)=Pr(N=k)%
\mathds{P}^k(\,\overset{\infty}{\underset{i=1}{\cup}}C_i)=Pr(N=k)\sum_{i=1}^%
\infty\mathds{P}^k(C_i)=\sum_{i=1}^\infty
Pr(N=k)\mathds{P}^k(C_i)=\sum_{i=1}^\infty\mathds{P}_V(A_i).
\end{equation*}
Therefore, by the $\pi$-$\lambda$ theorem \citep{Billingsley1995}, \ $\mathds{P}_V$
has a unique extension on $\sigma(\Pi)$.
\end{remark}

For designed experiments, we also may have a fixed sequence of design variables
$(X_{ij})_{j=1}^\infty$ so that the cluster variables are
$$(N_i; Y_{i1},Y_{i2},\ldots; X_{i1},X_{i2},\ldots )$$
but only $N_i$, $Y_i=(Y_{i1},\ldots,Y_{iN_i})'$
and $X_i=(X_{i1},\ldots,X_{iN_i})'$ are observed. Inference on the conditional
distribution of $Y_{ij}|X_{ij}$ may again be confounded by the cluster size.

\subsection{Parametric models}

Informative cluster size settings frequently appear in the biomedical
literature. There are  three natural ways to generate these type of
data. If the parametric model for the data generating mechanism can be correctly identified,
maximum likelihood estimates and likelihood ratio tests can be employed for statistical inference.

\begin{enumerate}
\item Models where the \emph{cluster size is assumed to have an influence on the outcomes}.
The data can then be thought to be generated in the following way:
First, $N_i$ is generated from its marginal distribution. Second,
$(Y_{ij})_{j=1}^\infty$ are generated from a conditional distribution
conditioned on $N_i$. If  the unobserved are
integrated out, the likelihood for what we observe is
$$ \prod_{i=1}^M P(N_i)  f(Y_{i1},\ldots,Y_{iN_i}|N_i). $$
See Remark 1.
These types of models are frequently found in the literature. The joint density
$f(Y_{i1},\ldots,Y_{iN_i}|N_i)$ can be for instance a multivariate normal with a intracluster correlation
coefficient $\rho(N_i)$.
\begin{example}[Fetal weights of mice, \citet{Dunson2003}]
\label{example:mousefetus} A female mouse is mated with a healthy
male likely to result in growing fetuses. If fewer fetuses are produced, more space and
nutritional resources will be available for those fetuses.
Therefore, there will be an inverse association between litter size
and the fetal weights.
\end{example}
\item Models where the  \emph{cluster size is assumed to depend on the outcomes}.
We suppose that the sequence $(Y_{ij})_{j=1}^\infty$
is first generated from its marginal distribution, and $N_i$ is then
generated from its conditional distribution conditioned on
$(Y_{ij})_{j=1}^\infty$. If $Y_{i1},Y_{i2},\ldots$ are i.i.d. and if
$P(N_i=k\vert Y_{i1},Y_{i2},\ldots)=P(N_i=k\vert
Y_{i1},\ldots,Y_{ik})$, $k=1,2,...$, the likelihood for what we
observe is
$$ \prod_{i=1}^M f(Y_{i1}, \ldots, Y_{iN_i}) P(N_i| Y_{i1},\ldots,Y_{iN_i}).$$
Note, however, that although $(Y_{ij})_{j=1}^\infty$ are
exchangeable, the observed variables $Y_{i1},\ldots,Y_{iN_i}$ may not
have this property.
The following example illustrates this kind of setting.
 \begin{example}
In the analysis of recurrent events during follow-up periods of fixed
lengths $c_i$, individuals (clusters) with a tendency to shorter
gaps contribute more events to the analysis than individuals
with a tendency to longer gaps. Now
$$ N_i=k \ \ \Leftrightarrow \ \ Y_{i1}+...+Y_{i,k-1}<c_i\le Y_{i1}+...+Y_{ik},$$
and the observation $Y_{ik}$ is right-censored.  As a result of the design,
Kaplan-Meier estimates are biased estimates of the
marginal survival function. Exchangeability condition
holds for $Y_{i1},...,Y_{i,k-1}$. Such designs are also subject to other
complexities; \citet{Huang2003} give a discussion.
\end{example}
\item In \emph{latent variable models}, a third unobservable random variable $\xi_i$
is assumed to be simultaneously influencing both the cluster size $N_i$ and the outcomes
 $Y_{i1},Y_{i2},\ldots$. The observed data likelihood contribution of the $i$th cluster is
obtained by integrating the latent variable $\xi_i$ and
$Y_{i,N_i+1},Y_{i,N_i+2},...$ out of the full likelihood expression. If $N_i$ and $Y_{i1},Y_{i2},\ldots$ are conditionally
independent and $Y_{i1},Y_{i2},\ldots$ conditionally i.i.d, then we
get the likelihood
$$\prod_{i=1}^M \left[\int p(N_i;\xi_i) \prod_{j=1}^{N_i} f(Y_{ij};\xi_i)  d Q(\xi_i)\right].$$
This is a model where exchangeability is met as well.

\begin{example}[\citet{WilliamsonDattaSatten2003}]
 Consider the association between the disease status of teeth, from a sample of individuals,
  and explanatory variables of interest. Disease status of teeth from the same person are correlated, and the individuals with poor dental status are likely to have fewer teeth. Poor dental health or hygiene can be thought to be the latent variable influencing both the disease status and the number of teeth.
\end{example}
\end{enumerate}

Note that introduction of treatments on the female mice in Example \ref{example:mousefetus}
may conceptually change it to a latent variable model. As nicely explained by \citet{Dunson2003},
treatments may have an effect on fetal weight with or without an effect on fetal losses.
Treatment may be acting on some unidentified latent variable, which influences both fetal weight
and fetal losses. They argue that in settings like this, it is important to model the cluster
size and the outcomes jointly, although the probability distribution of the cluster size is rarely of direct interest.

\subsection{Nonparametric models}

We have seen in the previous subsection that under
appropriate assumptions parametric models can be used. In can be argued,
though, that these conditions along with the usual distributional assumptions
are fairly restrictive, and perhaps unrealistic. In some settings it may be
difficult to postulate a model with natural parameters.

In the nonparametric approach, the general aim is to
make inference on the distributions of $Y_{i1},\ldots,Y_{iN_i}$ with an unspecified data generating
mechanism.  \citet{Bickel1975a,Bickel1975b}
introduced the general idea that one should first define measures (functionals)
of different interesting characteristics of the population and then use the
corresponding sample statistics as estimators. For example, if
$Y_1,...,Y_M$ is a random sample from an unknown univariate
distribution $F$, then for the sample mean $\bar Y=\frac 1 M
\sum_{i=1}^M Y_i$ and the sample variance $S^2=\frac 1M \sum_{i=1}^M
(Y_i-\bar Y)^2$, the corresponding functionals are the mean
functional $\mu(F)=\int xdF(x)$ and the variance functional
$\sigma^2(F)=\int(x-\mu(F))^2dF(x)$, respectively. Under general
assumptions on $F$, two functionals $\mu_1(F)$ and $\mu_2(F)$ may
recover the same value (\emph{e.g.} the mean functional and the median functional in
the nonparametric model of symmetric distribution) but the
statistical properties (\emph{e.g.} efficiency and robustness)  of the corresponding estimates (the sample
mean and the sample median) may be completely different.
The functional approach has now been generally adapted in the robust community.
  The approach seems particularly attractive in the framework of informative cluster size as it avoids much
of the inevitable complexities in the parametric approach.

\section{Quantities of interest as functionals}

\label{sec:functionals}

So far, the literature around informative cluster size problems has focused mainly on
regression problems and less on the inference on the marginal distribution of the outcome.
This section fills the apparent gap and defines a range of useful
marginal quantities for informative cluster size settings.

\subsection{Mean and variance functionals}

Let $V_1,\ldots,V_M$ be a random sample from distribution $\mathds{P}(V)$.
Commonly, the population functional of interest is a marginal
expected value. Under exchangeability, or if $Y_{i1},Y_{i2},\ldots$ are identically distributed
conditionally on $N_i$, the target parameter is then $E(Y_{i1})$. However, this assumption
is not needed throughout, and is relaxed in Remark 2.

%Of course, as $V_i$ are not fully observed, the functional should be

%defined so that the corresponding sample statistic $\hat

%T=T(\mathds{P}_M)$ can be observed.

A popular approach to estimate the parameter of interest has been to sample
one observation from each cluster randomly, apply standard methods for these i.i.d.
observations,  and ``average over'' all resampled sets \citep{Hoffman2001}.
Alternative resampling strategies have also been proposed \citep{Chiang2008}.
Indirectly but essentially this approach weights observations from different clusters with the inverses of the cluster size.
Another way of doing this would be to take the first observation of each cluster. Or to take all of them, and assign weights.

Our approach first identifies functionals, which recover the value $E (Y_{i1})$,
and secondly, chooses among the corresponding sample statistics. Consider functionals
$T(\mathds{P})$ for the unknown probability measure $\mathds{P}$.
Assuming that $Y_{i1},...,Y_{iN_i}$ are exchangeable, possible functionals are,
for example,
$$
T_1(\mathds{P})=E (Y_{i1})\ \ \mbox{and} \ \ T_2(\mathds{P})=
E\left[ \frac 1{N_i} \sum_{j=1}^{N_i} Y_{ij}\right],
$$
which lead to sample statistics
$$
\hat T_1=\frac 1M \sum_{i=1}^M Y_{i1}\ \ \mbox{and} \ \ \hat
T_2=\frac 1M \sum_{i=1}^M \frac 1{N_i} \sum_{j=1}^{N_i} Y_{ij},
$$
respectively. Although $T_1$ and $T_2$ are equal at the population level under the
exchangeability or conditional identical distribution assumption,
the corresponding sample statistics are very different. Note that under general assumptions,
$$\sqrt{M} \left ( \hat T_2 - T_2(\mathds{P}) \right ) \rightarrow_D N \left(0,\tau^2(\mathds{P}) \right ),$$
where
\begin{eqnarray*}
\tau^2(\mathds{P}) & = & E \left [ \left (\frac{1}{N_i} \sum_{j=1}^{N_i} Y_{ij} \right ) ^2 \right ] - T_2(\mathds{P})^2.
\end{eqnarray*}
A consistent estimate of $\tau^2(\mathds{P})$ is
\begin{eqnarray*}
\hat \tau^2 & = & \frac{1}{M} \sum_{i=1}^M \left (\frac{1}{N_i} \sum_{j=1}^{N_i} Y_{ij} \right ) ^2  - \hat T_2^2.
\end{eqnarray*}

\begin{figure}
\includegraphics[width=\textwidth]{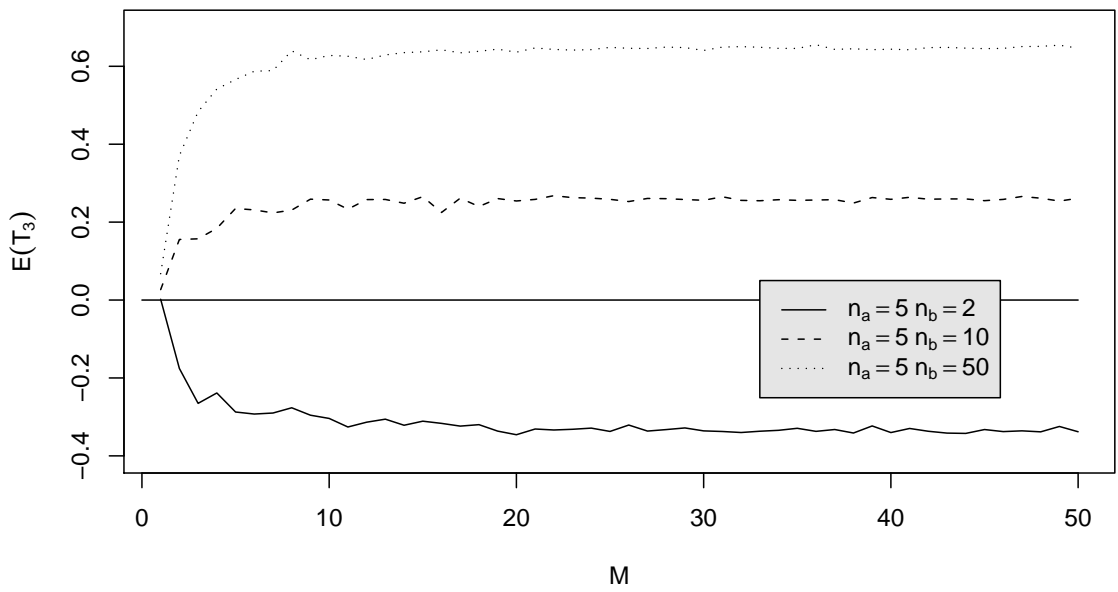}
\caption{\label{figure:samplemean} Expected value of the regular sample mean as a function of $M$ in the setting
of Example \ref{example:samplemean}. A numerical estimate of the expected value is shown.}
\end{figure}

\begin{example}[Population mean]
\label{example:samplemean} Consider the model $$Y_{ij} = \mu_{i} +
\epsilon_{ij},$$ where $\mu_i \sim N(0,1)$ and $\epsilon_{ij} \sim
N(0,1)$ independently, and the cluster sizes are generated via $N_i
= I(\mu_i < 0) n_a + I(\mu_i \geq 0) n_b$. The dependency of cluster
size on the realizations of $\mu_i$ induces the informative cluster
size property. We wish to estimate the population mean
$E(Y_{i1})=E(Y_{ij})$ which is zero in this case. $\hat T_1$ and
$\hat T_2$ are naturally unbiased, that is,
$$
E \left ( \frac{1}{M} \sum_{i=1}^M Y_{i1} \right )  = E \left (
\frac{1}{M} \sum_{i=1}^M \frac{1}{N_i} \sum_{j=1}^{N_i} Y_{ij}
\right ) = 0.
$$
The regular sample mean
$$
\hat T_3= \frac{1}{N} \sum_{i=1}^M \sum_{j=1}^{N_i} Y_{ij},
$$
where $N=\sum_{i=1}^M N_i$, is generally not unbiased and not even
consistent for $E(Y_{i1})$. In fact, there is
no corresponding functional $T_3(\mathds{P})$. The expected value
of $\hat T_3$ depends not only on $n_a$ and $n_b$, but also on $M$
(see Figure \ref{figure:samplemean}).

The variances of the two unbiased
estimators $\hat T_1$ and $\hat T_2$ are, of course, quite different:
$$Var\left ( \frac{1}{M} \sum_{i=1}^M Y_{i1} \right ) = 0.100 \hbox{ and }
Var\left ( \frac{1}{M} \sum_{i=1}^M \frac{1}{N_i}
\sum_{j=1}^{N_i} Y_{ij} \right ) \approx 0.057 ,$$
and we conclude that
$$\hat T_2 = \frac{1}{M} \sum_{i=1}^M \frac{1}{N_i} \sum_{j=1}^{N_i}
Y_{ij}$$ is the preferred estimate for the population quantity $E
(Y_{i1})$.
\end{example}

\begin{remark}
If it cannot be assumed that $Y_{i1},...,Y_{iN_i}$ are exchangeable,
a location center can still be defined as an expected value of a
randomly chosen observation in a random cluster, which still serves as a descriptive
statistic of the distribution. This functional is
then defined in the sample space of
$$V_i^* = (N_i; Y_{i1},\ldots,Y_{iN_i} ;\alpha_i)$$
where $\alpha_i$ is a pseudo random variable uniformly distributed
in $\{1,...,N_i\}$. Also $\alpha_i$ and $Y_{i1},\ldots,Y_{iN_i}$ are
conditionally independent conditioned on $N_i$. Let $\mathds{P}^*$
be the resulting probability measure. Then the location functional
is $T(\mathds{P}^*)=E^*(Y_{i\alpha_i})$. It is then straightforward
to see that
$$E^*\left[Y_{i\alpha_i}\right]=E^*\left[E^*(Y_{i\alpha_i}|N_i)\right]
 =E\left[\frac
1{N_i}\sum_{j=1}^{N_i} Y_{ij}\right].
$$
If $Y_{i1},...,Y_{iN_i}$ are exchangeable, then $T(\mathds{P}^*)=T_2(\mathds{P})$ is naturally
equal to $E(Y_{i1})$. The resulting estimate
is again
$$\frac 1M \sum_{i=1}^M \frac 1{N_i} \sum_{j=1}^{N_i} Y_{ij}.$$
This can be generalized into cases where $\alpha_i$ is not necessarily uniformly distributed.
\end{remark}

\begin{remark}
If the cluster size is noninformative, $Y_{i1},...,Y_{iN_i}$
exchangeable, and we wish to estimate $E(Y_{i1})$, all the weighted
means
$$\left (\sum_{i=1}^M\sum_{j=1}^{N_i} w_{ij} \right )^{-1} \sum_{i=1}^M\sum_{j=1}^{N_i}
w_{ij} Y_{ij}
$$
are unbiased. Optimal weights can be found in some simple settings.
Recall that, for informative cluster size, weights proportional
to $N_i^{-1}$ in the $i$th cluster guarantee unbiasedness. It may be possible to find
 other classes of weights that would also result in estimates having this property.
\end{remark}

In a similar way, there are several possible functionals for the
variance of $Y_{i1}$. Again,
$$\hat S_3= \frac{1}{N} \sum_{i=1}^M \sum_{j=1}^{N_i} (Y_{ij}-\hat
T_3)^2
$$
is biased,
$$
\hat S_1=\frac 1M \sum_{i=1}^M (Y_{i1}-\hat T_1)^2
$$
is appropriate under exchangeability but loses information, and
% $$S_2(\mathds{P})=\frac 1{N_i}

%\sum_{j=1}^{N_i} E\left[ (Y_{ij}-T_2(\mathds{P}))^2\right]$$
$$\hat S_2 =\frac 1M \sum_{i=1}^M  \frac 1{N_i}
\sum_{j=1}^{N_i}  (Y_{ij}-\hat T_2)^2$$ is the most natural
modification of the sample  variance  in this setting. It is only asymptotically unbiased. Unlike
in the standard setting, a general correction term to correct for the bias in finite samples cannot
be given. The variance functional corresponding to $\hat S_2$ is naturally
$$S_2(\mathds{P})= E\left[  \frac 1{N_i}
\sum_{j=1}^{N_i} (Y_{ij}-T_2(\mathds{P}))^2\right].$$

\subsection{Other distribution functionals}

Let $F$ denote the marginal distribution of a randomly chosen $Y$ in a random cluster. The
cumulative distribution function $F(y)$ and the corresponding
quantiles $q_{\alpha}=F^{-1}(\alpha)$ still serve as useful summary
measures of the marginal distribution even when exchangeability assumption is violated (Remark 2).
A natural functional for $F(y)$
is
$$
 F(y)(\mathds{P})=E \left [\frac{1}{n_{i}} \sum_{j=1}^{n_{i}} I(Y_{ij} \le y) \right
 ].
$$
The corresponding quantile functional $q_{\alpha}(\mathds{P})$
satisfies
$$
q_{\alpha}=\inf\{y\ : F(y)\ge \alpha  \}.
$$
The properly estimated cumulative distribution function yield
corresponding estimates of quantiles as well. The conventional estimate
$\frac{1}{N} \sum_{i=1}^M \sum_{j=1}^{N_i} I(Y_{ij}\leq x)$
is biased for these purposes and there is no corresponding functional.

Even though the correct functional form may appear obvious, we stress that
this simple functional structure must be maintained in all levels
when the functional nests other functionals within it. This need for
caution can be demonstrated by investigating the correct functional
of the $\alpha$-trimmed mean
$$E \left [ \frac{1}{(1-\alpha)N_i} \sum_{j=1}^{N_i} I(q_{\alpha/2} \leq Y_{ij} \leq q_{1-\alpha/2}) Y_{ij} \right ],$$
where the $q_{\alpha/2}$ and $q_{1-\alpha/2}$ are the corresponding
quantiles.  Now the quantiles themselves are functionals which
should be based on the correctly defined cumulative distribution functionals.

Another example is the sample
correlation, where the correct functional form employs three
redefined functionals, as shown by Example
\ref{example:correlation}.

\begin{example}[Sample correlation coefficient]
\label{example:correlation} Suppose that we observe a sample of i.i.d. clusters with bivariate observations
$$\{N_i; Y_{i1},\ldots, Y_{iN_i}\}, i = 1,\ldots,M,$$ where
$Y_{ij} = \mu_{i} + \epsilon_{ij}$ with $\mu_{i} \sim
N_2({0},{I}_2)$ and $\epsilon_{ij} \sim N_2({0},{I}_2)$. Thus,
$Y_{ij1}$ and $Y_{ij2}$ are uncorrelated $(Cov(Y_{ij1},Y_{ij2})=0)$
but the cluster size is informative in the following way. Assume
that large cluster sizes appear only when both components are large,
$$N_i = n_a + I(\mu_{i1}>1)I(\mu_{i2}>1)(n_b-n_a),$$
where $n_a=1$ and $n_b=10$. Consider the
biases of different sample statistics for $Cov(Y_{ij1},Y_{ij2})$ at $M=100$.
First,
\begin{eqnarray*}
E \left \{ \frac{1}{N} \sum_{i=1}^{M} \sum_{j=1}^{N_i} \left [
(Y_{ij1}-\bar{Y}_{1}) (Y_{ij2}-\bar{Y}_{2}) \right ] \right \} & = &
0.30,
\end{eqnarray*}
where $$
\bar{Y} = (\bar{Y}_{1},\bar{Y}_{2})' = \frac 1N \sum_{i=1}^M \sum_{j=1}^{N_i}
Y_{ij}.
$$
Deviation from zero indicates linear dependency and thus, it is clearly not a good
estimate of the marginal covariance functional of interest. The
weighted covariance functional
\begin{eqnarray*}
E \left \{\frac{1}{N} \sum_{i=1}^M \sum_{k=1}^{N_i}
\left [ (Y_{ij1}-\bar{Y}_{1}) (Y_{ij2}-\bar{Y}_{2}) \right ] \right
\} & = & 0.08,
\end{eqnarray*}
is still off the target because of the biased location estimate, and the correct covariance functional,
\begin{eqnarray*}
E \left \{\frac{1}{M} \sum_{i=1}^M \frac{1}{N_i} \sum_{k=1}^{N_i}
\left [ (Y_{ij1}-\tilde{Y}_{1}) (Y_{ij2}-\tilde{Y}_{2}) \right ]
\right \} & = & 0.00
\end{eqnarray*}
where the estimate of the mean vector is
$$
\tilde Y = (\tilde{Y}_{1},\tilde{Y}_{2})' = \frac 1M \sum_{i=1}^M \frac 1{N_i} \sum_{j=1}^{N_i}
Y_{ij}.
$$
The appropriate covariance functional is
$$E \left \{ \frac{1}{N_i} \sum_{j=1}^{N_i} \left [ Y_{ij1}-E \left (\frac{1}{N_i} \sum_{j=1}^{N_i} Y_{ij1} \right) \right ] \left [
Y_{ij2}-E \left (\frac{1}{N_i} \sum_{j=1}^{N_i} Y_{ij2} \right)\right ] \right \}.$$
To estimate the unknown correlation coefficient $Corr(Y_{ij1},Y_{ij2})$,
the correct covariance functional should be standardized using the square root
of the product of appropriately defined marginal variance functionals.
\end{example}

Example \ref{example:correlation} shows that incorrect functional forms can indicate dependency when it is actually
an artefact caused by informative cluster size. The opposite, incorrect functionals suggesting no dependency in presence
of real correlation, could also happen in a setting where the cluster sizes are related to the outcomes in a specific way.

\subsection{Sign and rank based functionals}

Much of the literature on informative cluster size has been on nonparametric methods. It
thus makes sense to define the basic sign and rank concepts and related quantities in the
same functional spirit.

Suppose that the interest lies in the median of the distribution
of $Y_{i1}$ rather than its expected value. Alternative location
functionals for this marginal median are then $T_1(\mathds{P})$ and
$T_2(\mathds{P})$ which satisfy
$$E \left [ I(Y_{i1} \leq T_1 ) \right ]=\frac 12
\ \ \mbox{and}\ \ E \left [ \frac{1}{N_i} \sum_{i=1}^{N_i} I(Y_{ij}
\leq T_2 ) \right ]=\frac 12,
$$
respectively. Then $\hat T_1$ is the sample median of $M$
observations  $Y_{11},...,Y_{M1}$ and $\hat T_2$ is the weighted
median of all $N=\sum_{i=1}^M N_i$ observations
$Y_{11},...,Y_{1N_1},Y_{21},..., Y_{MN_M}$ with weights proportional
to $1/N_i$ in the $i$th cluster.
 In the sign test one confronts the hypotheses
$$H_0:  T_2(\mathds{P}) = 0 \hbox{ vs. } H_1: T_2(\mathds{P})\neq 0.$$
A modified sign test statistic related to $T_2$ is
$$\frac{1}{M} \sum_{i=1}^M \frac{1}{N_i} \sum_{j=1}^{N_i} \sign(Y_{ij}),$$
which is a special case of  weighted sign tests considered by
\citet{LarocqueNevalainenOja2007}, but in the case of noninformative
cluster size, with weights chosen as $1/N_i$. Again, standard
asymptotics show that the standardized and quadratic form of the
test statistic has a limiting chi-square distribution with one
degree of freedom under the null hypothesis.

A signed-rank test for informative cluster size problems was
considered by \citet{Datta2008}. Their test is based on the
within-cluster resampling approach proposed by \citet{Hoffman2001},
which cleverly avoids the modeling of the covariance structure.
The concepts of rank and signed-rank can
be defined in an informative cluster size setting in a functional
manner. The estimate for the cumulative distribution functional $F(y)(\mathds{P})$ is
$$
 \hat F(y)= \frac 1M \sum_{i=1}^M \left [\frac{1}{N_{i}} \sum_{j=1}^{N_{i}} I(Y_{ij} \le y) \right ]
$$

An apparent modification of the  rank of $Y_{ij}$ is, accordingly, $\hat F(Y_{ij})$. For  the signed-rank concept we first need
$$
 F^+(y)(\mathds{P})=E \left [\frac{1}{N_{i}} \sum_{j=1}^{N_{i}} I(|Y_{ij}| \le y) \right ]
$$
and then the signed-rank of $Y_{ij}$ is $\sign(Y_{ij})\hat
F^+(|Y_{ij}|)$.   The signed-rank test statistic for testing whether
the symmetry center of the distribution of $Y_{ij}$ is zero is then
 $$\frac{1}{M} \sum_{i=1}^M \frac{1}{N_i} \sum_{j=1}^{N_i} \sign(Y_{ij}) \hat F^+ (|Y_{ij}|). $$
This test has been proposed in \citet{Datta2008}, where also the
limiting distribution is found via the resampling strategy.

To define the Hodges-Lehmann location functional for the
distribution of $Y_{i1}$ we need, as in the i.i.d. case, two
independent copies of the distribution of $V$, say, $V_i$ and
$V_{i'}$. Alternative location functionals $T_1(\mathds{P})$ and
$T_2(\mathds{P})$ satisfy
$$E \left [ I(Y_{i1}+Y_{i'1} \leq 2 T_1 ) \right ]=\frac 12
\ \ \mbox{and}\ \ E \left [ \frac{1}{N_i} \frac 1{N_{i'}}
\sum_{i=1}^{N_i} \sum_{j'=1}^{N_{i'}} I(Y_{ij}+Y_{i'j'} \leq  2 T_2
) \right ]=\frac 12,
$$
respectively. Note then that $\hat T_1$ is the Hodges-Lehmann
estimate calculated from  $M$ observations $Y_{11},...,Y_{M1}$ and
$\hat T_2$ is the weighted Hodges-Lehmann estimate based on  all $N$
observations. The signed-rank test given above is related to the
latter functional.

\subsection{Regression $L_2$ functionals}

In the linear regression case the cluster variables are
$(N_i; Y_{i1},Y_{i2},\ldots; X_{i1},X_{i2},\ldots )$.
On each cluster, we observe cluster sizes $N_i$, a matrix $Y_i=(Y_{i1},\ldots,Y_{iN_i})'$ consisting of multivariate
outcomes and a matrix of covariates $X_i=(X_{i1},\ldots,X_{iN_i})'$. Assume that $$(X_{i1},Y_{i1}),...,
(X_{iN_i},Y_{iN_i})$$ are exchangeable and  that, for the correct
$\beta$, $E\left((Y_{i1}-\beta'X_{i1})X_{i1}' \right)=0$. Possible regression coefficient functionals in this case are
$\beta_1(\mathds{P})$ and $\beta_2(\mathds{P})$ satisfying
$$
E\left((Y_{i1}-\beta_1'X_{i1})X_{i1}' \right)=0 \ \ \mbox{and}\ \
E\left(\frac 1{N_i}\sum_{j=1}^{N_i} (Y_{ij}-\beta_2'X_{ij})X_{ij}'
\right)=0.
$$
Then $\hat\beta_1$ is the least squares estimate based on
$(X_{11},Y_{11}),..., (X_{M1},Y_{M1})$ and $\hat\beta_2$ is a
weighted least squares estimate using all the observations in the appropriate way.
Again, in the case of informative sample size, there is no functional
corresponding to the naive estimate $\hat\beta$ satisfying
$$ \sum_{i=1}^M \sum_{j=1}^{N_i} (Y_{ij}-\hat\beta'X_{ij})X_{ij}'=0. $$
It is useful to note that from the estimating equation
$$\frac{1}{M}\sum_{i=1}^M \frac{1}{N_i} \sum_{j=1}^{N_i}(Y_{ij}-\hat\beta_2' X_{ij})X_{ij}' = 0$$
we get by straightforward calculation that
$$\hat\beta_2 = \left [ \frac{1}{M} \sum_{i=1}^{M} \frac{1}{N_i} X_i'X_i\right ]^{-1}
\left [ \frac{1}{M} \sum_{i=1}^{M} \frac{1}{N_i} X_i'Y_i\right ],$$
and that
$$\sqrt{M}\left ( \hat\beta_2 - \beta_2 \right ) = \left [ \frac{1}{M} \sum_{i=1}^{M} \frac{1}{N_i} X_i'X_i\right ]^{-1} \left [ \frac{1}{\sqrt{M}} \sum_{i=1}^{M} \frac{1}{N_i} X_i'R_i\right ],$$
where $R_i=(R_{i1},...,R_{iN_i})'$ with $R_{ij}=Y_{ij}-\beta_2'X_{ij}$, $j=1,...,N_i$.
The first part of the right hand side converges in probability to its expected value, and the second part clearly has
a limiting normal distribution. Thus, the statistical inference on the regression coefficients
$( \hat\beta_2 - \beta_2 ) $ can be based on a
normal distribution with mean zero and a covariance matrix, which can be estimated by the sandwiching form
$\hat A^{-1}\hat B \hat A^{-1}$, where
$$
\hat A  =  \sum_{i=1}^M \frac{1}{N_i}X_{i}'X_{i} \hbox{\; and \; }
\hat B  =  \sum_{i=1}^M \left ( \frac{1}{N_i} X_{i}'R_i \right ) \left ( \frac{1}{N_i} X_{i}'R_i \right )'.
$$

For cases where $(X_{i1},Y_{i1}),...,
(X_{iN_i},Y_{iN_i})$ are not exchangeable, we can define
$\beta(\mathds{P}^*)$ in the sample space of
$$V_i^* = (N_i; Y_{i1},\ldots,Y_{iN_i}; X_{i1},\ldots,X_{iN_i}; \alpha_i ).$$
by
$$
E^*\left((Y_{i\alpha_i}-\beta' X_{i\alpha_i})X_{i\alpha_i} \right)=
E\left(\frac 1{N_i}\sum_{j=1}^{N_i} (Y_{ij}-\beta' X_{ij})X_{ij}
\right)=0.
$$
 Functional $\beta(\mathds{P}^*)$ then is a measure of ``average
 regression''.

In standard cluster specific random effects models, if the cluster size distribution only depend on the random effect but not on the covariates,
then a simple calculation shows that the corresponding components of the estimators without the $1/N_i$ weight also converge to the correct
regression parameters. In other words, in such cases, the classical analyses will also work for such parameters but the estimates of the other
parameters including the intercept terms will continue to be biased. \citet{Benhin2005}, \citet{Gueorguieva2005}, \citet{Wang2011} and
\citet{Neuhaus2011}, e.g., consider these types of results.

\subsection{Regression M functionals}

In the univariate case, simultaneous (naive) M functionals $\beta(\mathds{P})$ and $\sigma(\mathds{P})$ for linear regression  are given by
equations
\[
E(w_1(R_{i1}) R_{i1}X_{i1})=0 \ \ \mbox{and} \ \ E(w_2(R_{i1})R_{i1}^2)=E(w_3(R_{i1})),
\]
where now
\[
R_{ij}=R_{ij}(\beta,\sigma)=\frac {Y_{ij}-\beta'X_{ij}} \sigma,\ \ i=1,...,M;\ j=1,...,N_i.
\]
Recall that, in maximum likelihood estimation,
\[
w_1(R)=w_2(R)= -\frac {f'(R)} {R f(R)} \ \ \mbox{and}\ \ w_3(R)\equiv 1,
\]
and Huber's estimate, for example, is given by
\[
w_1(R)=\min\left(1,c/|R| \right), \ \
w_2(R)=d\min\left(1,c^2/R^2 \right),\ \ \mbox{and}\ \
w_3(R) \equiv1,
\]
with tuning parameters $c,d>0$. Alternative functionals for clustered data with informative cluster size are given by
\[
E\left( \frac 1{N_i}\sum_{j=1}^{N_i}w_1(R_{ij}) R_{ij} X_{ij}\right)=0 \ \ \mbox{and} \ \ E\left(\frac 1 {N_i} \sum_{j=1}^{N_i} w_2(R_{ij} )R_{ij}^2\right)=E\left(\frac 1 {N_i} \sum_{j=1}^{N_i}w_3(R_{ij})\right),
\]
Iteration steps to compute the estimates $\hat\beta$ and $\hat\sigma$
are:
\begin{enumerate}
\item first update the residuals
\[
R_{ij} \leftarrow \frac {Y_{ij}-\hat{\beta}'X_{ij}} {\hat\sigma},
\ \  R_{i} \leftarrow(R_{i1},...,R_{iN_i})',
\]
\item next the weights
\begin{eqnarray*}
% \nonumber to remove numbering (before each equation)
  W_{1i} &\leftarrow & diag(w_1(R_{i1}),..., w_1(R_{iN_i}))   \\
  W_{2i} &\leftarrow & diag(w_2(R_{i1}),..., w_2(R_{iN_i}))   \\
  W_{3i} &\leftarrow & diag(w_3(R_{i1}),..., w_3(R_{iN_i}))
\end{eqnarray*}
\item and finally obtain new values of $\hat\beta$ and $\hat\sigma$ as
\begin{eqnarray*}
% \nonumber to remove numbering (before each equation)
  \hat \beta &\leftarrow & \left( \frac 1M \sum_{i=1}^M \frac 1{N_i} X_i'W_{1i} X_i \right)^{-1}
  \left( \frac 1M \sum_{i=1}^M \frac 1{N_i} X_i'W_{1i} Y_i \right)
     \\
  \hat\sigma^2 &\leftarrow &  \left( \frac 1M \sum_{i=1}^M \frac 1{N_i} 1_{N_i}'W_{3i} 1_{N_i} \right)^{-1}
  \left( \frac 1M \sum_{i=1}^M \frac 1{N_i} R_i'W_{2i} R_i \right)\hat\sigma^2.
\end{eqnarray*}
\end{enumerate}
The covariance matrix estimate of $\hat\beta$ can be approximated by the sandwich estimate
$\hat A^{-1}\hat B \hat A^{-1}$ where now
\[
\hat A=  \sum_{i=1}^M \frac 1{N_i} X_i'W_{1i} X_i
\ \ \mbox{and}\ \
\hat B=  \sum_{i=1}^M \left(\frac 1{N_i} X_i'W_{1i} R_i\right) \left(\frac 1{N_i} X_i'W_{1i} R_i\right)'
\]

\section{Test construction}

\label{sec:testing}

Suppose that the null hypothesis of interest $H_0$ implies that
(or can be formulated as)
\[
E^*[ T(Y_{1\alpha_1},...,Y_{M\alpha_M})]=0
\]
where $E^*$  corresponds to the probability measure $\mathds{P}^\ast$
 overlying the distribution of
$\alpha_i,\,N_i,Y_{i1},\ldots,Y_{iN_i}$, plus the covariates $X_{i1},\ldots,X_{iN_i}$
in the regression case, if they are assumed to be random, too. An appropriate test statistic for this testing problem
based on i.i.d. observations $ Y_{1\alpha_1},...,Y_{M\alpha_M}$ (plus covariates
in case of regression) is simply given by
\[
\hat{T}^*=  T(Y_{1\alpha_1},...,Y_{M\alpha_M}).
\]
As for example, if we are interested in testing $H_0:\ E^*(Y_{i\alpha_i})=0$, then
the null hypothesis implies that also $E^*[ T(Y_{1\alpha_1},...,Y_{M\alpha_M})]=0$ with
 $T(Y_{1\alpha_1},...,Y_{M\alpha_M})=M^{-1}\sum_i  Y_{i\alpha_i}$. While
$\hat{T}^*$ is a valid test statistic, it is objectionable: (i) this may be inefficient since
a large part of the data will be ignored
depending on which observations are chosen by the particular realization of
the random indices $\alpha_i$ and (ii) the artificial randomization
itself may be unsatisfactory for practical application and may lead to
additional variability. Therefore, an appropriate strategy will be to take a
further expectation of this test statistic $\hat{T}^*$ with respect to
conditional distribution of the indices $\alpha_i$ given the original
clustered data $V_1,\ldots, V_M$ leading to the test statistic
\begin{equation*}
\hat{T}(V_1,\ldots,V_M)=E^*(\hat{T}^*|V_1,\ldots,V_M).
\end{equation*}

Depending on the problem, this can sometimes be analytically calculated
exactly or approximately (up to terms that are asymptotically ignorable, as
$M \rightarrow\infty$)
through a linear approximation of $\hat{T}^*$ \citep{DattaSatten2005,Datta2008}.
In the one sample problem of
testing location symmetry, \citet{Datta2008}  adopted the signed-rank
statistic for clustered data by selecting $\hat{T}^*\,$ to be the
regular signed-rank statistic for i.i.d. data. It turns out that the
resulting $\hat{T}\,$is algebraically equivalent to the signed-clustered
rank test statistic we obtained in Section \ref{sec:functionals} from an intuitive consideration
via statistical functional.

The test statistic $\hat{T}$ can always be estimated using a Monte-Carlo
technique that is in the same spirit of the original within-cluster resampling proposal by \citet{Hoffman2001}
for the estimation problem:
\begin{eqnarray*}
\hat{T}\approx\frac{1}{B}\sum_{b=1}^B\hat{T}^*(Y_{1\alpha_1(b)},\ldots,Y_{M\alpha_M(b)}),
\end{eqnarray*}
where a large number $B$ sets of realizations of the random indices
$(\alpha_1,\cdots\,\alpha_M)\,$ are drawn.

An estimate of the sampling variance of this test statistic can be computed in one of three
possible ways:
\begin{itemize}
\item[(i)] by analytical calculation involving linearization
techniques such as projections;
\item[(ii)] by the Monte-Carlo variance formula
\begin{eqnarray*}
\widehat{Var}(\hat{T}) & \approx & \frac{1}{B}\sum_{b=1}^B \widehat{Var}
\left \{
\hat{T}^*(Y_{1\alpha_1(b)},\ldots,Y_{M\alpha_M(b)})\right \}\\
& & -\frac{1}{B-1}\sum_{b=1}^B \left \{ \hat{T}^*(Y_{1\alpha_1(b)},\ldots,\,Y_{M\alpha_M(b)})-\hat{T} \right \}^2,
\end{eqnarray*}
with the assumption that one has a variance formula for the statistic
$\hat{T}^*$; or
\item[(iii)] by bootstrap resampling of the entire cluster of
observations $V_i$ and by empirical variance of the test values of the test
statistics calculated with the resampled data.
\end{itemize}

\begin{example}[A modified \textit{t}-test]
An immediate modification of the $t$-test in conjunction with informative
cluster size is as follows. The goal is to confront the null
hypothesis $H_0: E^\ast(Y_{i\alpha_i})=0$ with the alternative $H_0: E^\ast(Y_{i\alpha_i})
\neq 0$, where $\alpha_i$ is uniformly distributed. Alternatively, the null can be formulated
as $H_0: E^\ast \left ( \frac{1}{M} \sum_{i=1}^M Y_{i\alpha_i} \right ) = 0$, which gives the test
statistic $E^\ast \left ( \frac{1}{M} \sum_{i=1}^M Y_{i\alpha_i} | V_1,\ldots, V_M \right ) = \hat T_2$.
As $M$ tends to infinity, the limiting distribution of the
modified one-sample $t$-statistic is
$$ {\sqrt{M} \hat T_2}/{\hat \sigma} \rightarrow_D N(0,1)$$
where
$$\hat\sigma^2=\frac{1}{M} \sum_{i=1}^M \frac{1}{N_i^2} \left (   \sum_{j=1}^{N_i} Y_{ij} \right)^2$$
is a consistent estimate of the limiting variance of $\sqrt{M} \hat
T_2$, because $g(V_i) = \frac{1}{N_i} \sum_{j=1}^{N_i} Y_{ij}$ are i.i.d. with expectation zero under the null.
\end{example}

\section{A data example on patients with incomplete spinal cord injuries}

\label{sec:dataexample}

\begin{figure}
\includegraphics[width=12cm]{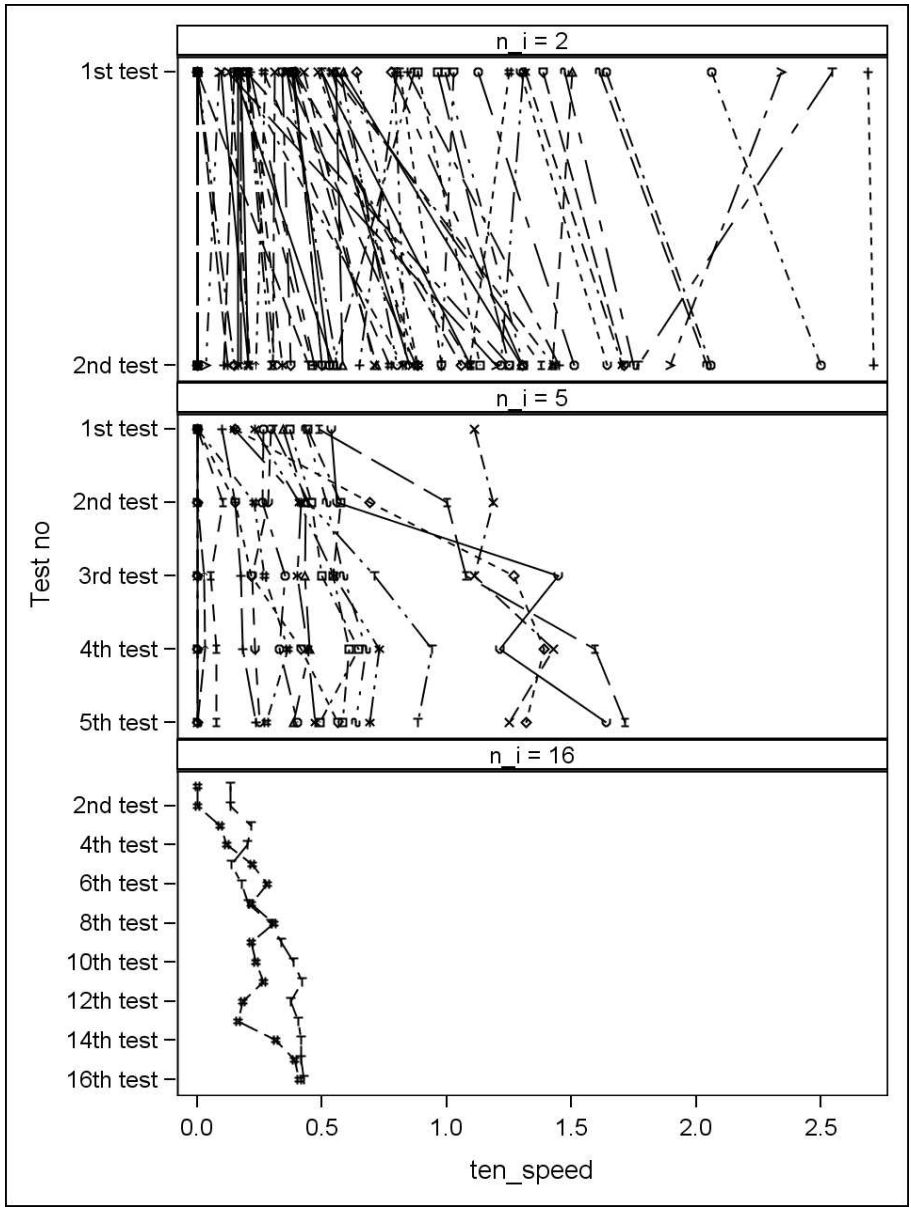}
\caption{\label{figure:scatter}Results in the Ten Meter Walk Test depend on the cluster size, the total number of tests. Results are shown for selected cluster sizes. Furthermore, the tendency for  improvement in tests results over time is clearly visible.}
\end{figure}

This data set is based on an observational cohort of patients at the NeuroRecovery Network (NRN). Patients eligible for NRN have incomplete spinal cord injuries (SCI) with lesion at level T10 or above and are not participating in inpatient rehabilitation programs \citep{Harkema2012}. Patients are discharged from the NRN for non-compliance with treatment, patient election, or if a plateau in the recovery of function is achieved.  This last discharge criterion is of particular interest to the present analysis and is the reason for the potential informativeness of the cluster size.  More severely impaired patients tended to have more ``room for improvement'' in function and hence remained enrolled in the NRN for longer periods of time, contributing more observations than those that enrolled with higher pre-existing function.  This phenomenon has been previously demonstrated for NRN patients (Figure \ref{figure:scatter}).

The outcome measures per longitudinal evaluation are as follows:
The Ten Meter Walk Test is commonly used as a measure of walking capacity in SCI patients.  In each test, a patient is instructed to walk as fast (10MW) as possible without assistance from the therapist conducting the assessment.  The reliability and validity of this test in measuring walking function has previously been demonstrated \citep{Hedel2005}.

The sample mean of the ten meter walking speed (meters/second) is $\hat T_3=0.439$. The weighted mean (weights inversely proportional
to cluster size) is $\hat T_2 = 0.493$. This indicates that the marginal mean is underestimated without the proper weighting;
in fact, $\hat T_3$ does not estimate any population functional and it is therefore not a proper estimate.
The sample statistic $\hat T_1=0.373$ does not estimate the same quantity as $\hat T_2$, either, because the assumption of
exchangeability is not reasonable due to improvement in patient conditions over time (Figure \ref{figure:scatter}).

We investigate for further illustration purposes the behavior of four estimates of $\beta$ in the
linear  model of the  form $$Y_{ij}=\beta'X_{ij}+\epsilon_{ij},$$
for studying the effects of gender (1 if male; 0 if female) and four races (indicators race$_1,\ldots$, race$_3$) on the results
of the Ten Meter Walk Test (the Ten Meter Walk speed).
There are $M=333$ patients with at least one test result
and $\sum_i N_i = 1329$ test results with non-missing values on the covariates.
The four estimates of regression coefficients and their standard errors are obtained as follows.
\begin{enumerate}
\item Within-cluster resampling (WCR) with 1000 resamples. The model is fitted using ordinary least squares on each resample.
Variances and standard errors of the regression coefficients are estimated by the Monte-Carlo variance formula given in section
\ref{sec:testing}.
\item Ordinary least squares (OLS) fitted on the full data completely ignoring the clustering.
\item Inverse cluster size weighted least squares (ICSWLS) with weights inversely proportional to cluster size, and fitted on the full data. Standard errors
are derived from the sandwich estimator of the covariance
 matrix given in section 3.
\item Linear mixed model (MM) with an additional patient random effect. Parameter estimates and their standard errors
are derived via restricted maximum likelihood estimation using inverse of the covariance matrix as the weight matrix.
\item Inverse cluster size weighted Huber's regression (ICWHR) estimate with $c=1.5$ and $d=1$ accompanied with standard errors derived from the corresponding
sandwich estimator of the covariance.
\end{enumerate}

\begin{table*}
\begin{footnotesize}
\caption{Parameter estimates and their standard errors (in parentheses) obtained by the
within-cluster resampling (WCR), ordinary least squares (OLS),
weighted least squares (ICSWLS), linear mixed model (MM) and weighted Huber's regression estimate (ICWHR).}
\begin{tabular}{crrrrrrrrrr}
\hline
& \multicolumn{10}{c}{Estimate (standard error)}\\
\cline{2-11}
Parameter        & \multicolumn{2}{c}{WCR} & \multicolumn{2}{c}{OLS} & \multicolumn{2}{c}{ICSWLS} & \multicolumn{2}{c}{MM} & \multicolumn{2}{c}{ICWHR}\\
\hline
Intercept        &  0.418 & (0.216) &  0.297 & (0.104) &  0.420  & (0.203)  &  0.415 & (0.211) & 0.385 & (0.169)\\
gender           &  0.139 & (0.064) &  0.107 & (0.029) &  0.138  & (0.059)  &  0.137 & (0.064) & 0.094 & (0.047)\\
race$_1$         &  0.086 & (0.220) &  0.172 & (0.107) &  0.086  & (0.218)  &  0.088 & (0.214) & 0.031 & (0.176)\\
race$_2$         &  0.042 & (0.244) &  0.258 & (0.119) &  0.040  & (0.229)  &  0.053 & (0.238) & 0.032 & (0.216)\\
race$_3$         & -0.060 & (0.211) &  0.038 & (0.102) & -0.061 &  (0.202)  & -0.056 & (0.207) & -0.063 & (0.167)\\
\hline
\end{tabular}
\end{footnotesize}
\end{table*}

WCR and ICSWLS approaches result into nearly identical parameter estimates and they are both known
to be unbiased. Differences are attributable to randomness arising from resampling. Their standard
errors are  similar throughout. OLS estimates are biased, and severely so.
The regression coefficient for race$_3$,  even has a different sign. Furthermore, the standard errors of the estimates are
artificial and much too small as they do not account for the clustering. The parameter estimates from the linear mixed model
often fall between ICSWLS and OLS estimates and are not far off, either. It has been noted that
under specific conditions, a linear mixed model can result
into consistent estimates of the slope parameters \citep{Neuhaus2011}; a finding that is supported
by these analyses. This, however, is not the case
even in a random intercepts model if the covariate is related to the cluster sizes \citep{Wang2011,Lorenz2011}.
In our setting the explanatory variables are not closely related to the cluster sizes and this is a
potential reason for the good performance of the linear mixed model here. Among these methods,
our preference would be the ICSWLS leading to unbiased estimates without computational burden due to resampling
of data. Inverse cluster size weighted Huber's estimate provides a robust alternative to these methods and performs extremely well
for this particular data set with similar estimates of regression coefficients and smaller standard errors
throughout.

\section{Concluding remarks}

\label{sec:concluding}

This paper gives an account on appropriate models, summary statistics and generalizing statistical
classical, nonparametric and robust procedures on clustered data with possibly informative cluster size. We have demonstrated
how subtle the problem is, and hope to have convinced the reader of a general method of dealing with it
by appropriate functionals, leading to weighted sample statistics.
In fact, it seems to us that the
whole classical statistical theory, and the theory of robust statistical
procedures with the concepts such as the breakdown point and influence
function can be reformulated along these lines for informative cluster
size  problems.

It is clear that
not all clustered data suffer from  informative cluster size. Nevertheless, it seems like
a good idea to investigate the distribution of the responses as function
of the cluster size by means of graphical summaries or similar, to make sure this is not
the case. Note, however, that the proposed modified properties are valid also if the cluster size is not informative,
at the possible cost of losing some efficiency relative to optimally weighted procedures.

\section*{Acknowledgements}

This research was supported by the Academy of Finland and by NIH grants 1R03DE020839-01A1, 5R03DE020839-02  and 1R03DE022538-01. 
The authors
are grateful for the use of data from the NeuroRecovery Network,
and thank the directors of centers participating in the NRN: Steve
Ahr (Frazier Rehab Institute, Louisville, KY), Steve Williams, MD
(Boston Medical Center, Boston, MA), Daniel Graves, PhD (Memorial
Hermann/The Institute of Rehabilitation and Research, Houston,
TX), Keith Tansey, MD, PhD (Shepherd Center, Atlanta, GA), Gail
Forrest, PhD (Kessler Medical Rehabilitation Research and
Education Corporation, West Orange, NJ), D. Michele Basso PT,
EdD (The Ohio State University Medical Center, Columbus, OH) and
Mary Schmidt Read, PT, DPT, MS (Magee Rehabilitation, Philadelphia, PA).

\bibliography{articles,books}

\begin{thebibliography}{21}
\providecommand{\natexlab}[1]{#1}
\providecommand{\url}[1]{\texttt{#1}}
\expandafter\ifx\csname urlstyle\endcsname\relax
  \providecommand{\doi}[1]{doi: #1}\else
  \providecommand{\doi}{doi: \begingroup \urlstyle{rm}\Url}\fi

\bibitem[Benhin et~al.(2005)Benhin, Rao, and Scott]{Benhin2005}
E.~Benhin, J.~N.~K. Rao, and A.~J. Scott.
\newblock Mean estimating equation approach to analysing cluster-correlated
  data with nonignorable cluster sizes.
\newblock \emph{Biometrika}, 92:\penalty0 435--450, 2005.

\bibitem[Bickel and Lehmann(1975{\natexlab{a}})]{Bickel1975a}
P.~J. Bickel and E.~L. Lehmann.
\newblock Descriptive statistics for nonparametric models. {I}. {I}ntroduction.
\newblock \emph{The Annals of Statistics}, 3:\penalty0 1038--1044,
  1975{\natexlab{a}}.

\bibitem[Bickel and Lehmann(1975{\natexlab{b}})]{Bickel1975b}
P.~J. Bickel and E.~L. Lehmann.
\newblock Descriptive statistics for nonparametric models. {II}. {L}ocation.
\newblock \emph{The Annals of Statistics}, 3:\penalty0 1045--1069,
  1975{\natexlab{b}}.

\bibitem[Billingsley(1995)]{Billingsley1995}
P.~Billingsley.
\newblock \emph{Probability and Measure}.
\newblock John Wiley \& Sons, New York, USA, third edition, 1995.

\bibitem[Chiang and Lee(2008)]{Chiang2008}
C.-T. Chiang and K.-Y. Lee.
\newblock Efficient estimation methods for informative cluster size data.
\newblock \emph{Statistica Sinica}, 18:\penalty0 121--133, 2008.

\bibitem[Datta and Satten(2005)]{DattaSatten2005}
S.~Datta and G.~A. Satten.
\newblock Rank-sum tests for clustered data.
\newblock \emph{Journal of the American Statistical Association}, 100:\penalty0
  908--915, 2005.

\bibitem[Datta and Satten(2008)]{Datta2008}
S.~Datta and G.~A. Satten.
\newblock A signed-rank test for clustered data.
\newblock \emph{Biometrics}, 64:\penalty0 501--507, 2008.

\bibitem[Dunson et~al.(2003)Dunson, Chen, and Harry]{Dunson2003}
D.~B. Dunson, Z.~Chen, and J.~Harry.
\newblock A {B}ayesian approach for joint modeling of cluster size and
  subunit-specific outcomes.
\newblock \emph{Biometrics}, 59:\penalty0 521--530, 2003.

\bibitem[Gueorguieva(2005)]{Gueorguieva2005}
R.~V. Gueorguieva.
\newblock Comments about joint modeling of cluster size and binary and
  continuous subunit-specific outcomes.
\newblock \emph{Biometrics}, 61:\penalty0 862--867, 2005.

\bibitem[Harkema et~al.(2012)Harkema, Schmidt-Read, Behrman, Bratta, Sisto, and
  Edgerton]{Harkema2012}
S.~J. Harkema, M.~Schmidt-Read, A.~Behrman, A.~Bratta, S.~A. Sisto, and V.~R.
  Edgerton.
\newblock Establishing the neurorecovery network: multi-site rehabilitation
  centers that provide activity based therapies and assessments for neurologic
  disorders.
\newblock \emph{Archives of Physical Medicine and Rehabilitation}, 93:\penalty0
  1498--1507, 2012.

\bibitem[Hoffman et~al.(2001)Hoffman, Sen, and Weinberg]{Hoffman2001}
E.~B. Hoffman, P.~K. Sen, and C.~R. Weinberg.
\newblock Within-cluster resampling.
\newblock \emph{Biometrika}, 88:\penalty0 1121--1134, 2001.

\bibitem[Huang and Chen(2003)]{Huang2003}
Y.~Huang and Y.~Q. Chen.
\newblock Marginal regression of gaps between recurrent events.
\newblock \emph{Lifetime Data Anal}, 9\penalty0 (3):\penalty0 293--303, Sep
  2003.

\bibitem[Larocque et~al.(2007)Larocque, Nevalainen, and
  Oja]{LarocqueNevalainenOja2007}
D.~Larocque, J.~Nevalainen, and H.~Oja.
\newblock A weighted multivariate sign test for cluster correlated data.
\newblock \emph{Biometrika}, 94:\penalty0 267--283, 2007.

\bibitem[Lorenz et~al.(2011)Lorenz, Datta, and Harkema]{Lorenz2011}
D.~J. Lorenz, S.~Datta, and S.~J. Harkema.
\newblock Marginal association measures for clustered data.
\newblock \emph{Statistics in Medicine}, 30\penalty0 (27):\penalty0 3181--3191,
  Nov 2011.
\newblock \doi{10.1002/sim.4368}.
\newblock URL \url{http://dx.doi.org/10.1002/sim.4368}.

\bibitem[Neuhaus and McCulloch(2011)]{Neuhaus2011}
J.~M. Neuhaus and C.~E. McCulloch.
\newblock Estimation of covariate effects in generalized linear mixed models
  with informative cluster sizes.
\newblock \emph{Biometrika}, 98:\penalty0 147--162, 2011.

\bibitem[Panageas et~al.(2007)Panageas, Schrag, Russell, Venkatraman, and
  Begg]{Panageas2007}
K.~S. Panageas, D.~Schrag, L.~A. Russell, E.~S. Venkatraman, and C.~B. Begg.
\newblock Properties of analysis methods that account for clustering in
  volume-outcome studies when the primary predictor is cluster size.
\newblock \emph{Statistics in Medicine}, 26:\penalty0 2017--2035, 2007.

\bibitem[van Hedel et~al.(2005)van Hedel, Wirz, and Dietz]{Hedel2005}
H.~van Hedel, M.~Wirz, and V.~Dietz.
\newblock Assessing walking ability in subjects with spinal cord injury:
  validity and reliability of 3 walking tests.
\newblock \emph{Archives of Physical Medicine and Rehabilitation}, 86:\penalty0
  190--196, 2005.

\bibitem[Wang et~al.(2011)Wang, Kong, and Datta]{Wang2011}
M.~Wang, M.~K. Kong, and S.~Datta.
\newblock Inference for marginal linear models for clustered longitudinal data
  with potentially informative cluster sizes.
\newblock \emph{Statistical Methods in Medical Research}, 20:\penalty0
  347--367, 2011.
\newblock \doi{doi: 10.1177/0962280209347043}.

\bibitem[Williamson et~al.(2003)Williamson, Datta, and
  Satten]{WilliamsonDattaSatten2003}
J.~M. Williamson, S.~Datta, and G.~A. Satten.
\newblock Marginal analyses of clustered data when cluster size is informative.
\newblock \emph{Biometrics}, 59:\penalty0 36--42, 2003.

\bibitem[Williamson et~al.(2007)Williamson, Kim, and Warner]{Williamson2007}
J.~M. Williamson, H.-Y. Kim, and L.~Warner.
\newblock Weighting condom use data to account for nonignorable cluster size.
\newblock \emph{Annals of Epidemiology}, 17:\penalty0 603--607, 2007.

\bibitem[Williamson et~al.(2008)Williamson, Kim, Manatunga, and
  Addiss]{Williamson2008}
J.~M. Williamson, H.-Y. Kim, A.~Manatunga, and D.~G. Addiss.
\newblock Modeling survival data with informative cluster size.
\newblock \emph{Statistics in Medicine}, 27:\penalty0 543--555, 2008.

\end{thebibliography}

\bibliographystyle{abbrvnat}

\newpage

\end{document}